# $\underline{\boldsymbol{Parametric\ solutions\ to\ (\ six)\ n^{th}\ powers\ equal\ to}}$

# $\underline{\boldsymbol{another\ (six)\ n^{th}\ powers\ \ for\ degree\ 'n'\ =\ 2,3,4,5,6,7,8,\&\ 9}}$

Authors: *Seiji Tomita*[1] *and Oliver Couto*[2]

1 Computer Engineer, Tokyo software co. Inc., Japan.
Email: fermat@m15.alpha-net.ne.jp

2 University of Waterloo, Waterloo, Ontario, Canada
Email: samson@celebrating-mathematics.com

Abstract

Consider the below mentioned equation:

$$[a^n + b^n + c^n + d^n + e^n + f^n] = [p^n + q^n + r^n + s^n + t^n + u^n] \quad \text{----- (A)}$$

Historically in math literature there are instances where solutions have been arrived at by different authors for equation (A) above. Ref.no. (1) by A. Bremner & J. Delorme and Ref. no. (10) by Tito Piezas. The difference is that this article has done systematic analysis of equation (A) for n=2,3,4,5,6,7,8 & 9. While numerical solutions for equation (A) is available on “Wolfram math” website, search for parametric solutions to equation (A) in various publications for all n=2,3,4,5,6,7,8 & 9 did not yield much success. The authors of this paper have selected six terms on each side of equation (A) since the difficulty of the problem increases every time a term is deleted on each side of equation (A). The authors have provided parametric solutions for equation (A) for n=2, 3, 4, 5 & 6 and for n=7, 8 & 9 solutions using elliptical curve theory has been provided. Also we would like to mention that solutions for n=7, 8 & 9 have infinite numerical solutions.

**Degree two, n=2**

$$a^2+b^2+c^2+d^2+e^2+f^2 = p^2+q^2+r^2+s^2+t^2+u^2 \text{ ----- (1)}$$

We have the below mentioned numerical solution:

$$(1,30,31,36,7,17)^2 = (3,4,19,27,34,35)^2$$

In the above let a=1, b=30 & c=31, Since 1+30=31 or (a+b=c)

Hence we have $(a^2+b^2+c^2) = 2(a^2+ab+b^2)$ -------------- (2)

Let a= (30+t) and b=(1+kt) .Substituting for (a,b) & solving for 't' in (a^2+ab+b^2) in

right hand side of equation (2) above , we get parameterization of (a & b) given below.

a= $(30k^2-2k-31)/(k^2+k+1)$

b= $(-31k^2-60k+1)/(k^2+k+1)$

c=$a+b = (-k^2-62k-30)/(k^2+k+1)$

Also since $(a,b,c)^2 = (p,q,r,s,t,u)^2 - (d,e,f)^2$ , after substituting above values of

(a,b,c) we get the below mentioned parameterization of equation (A) for degree two.

$$(30k^2-2k-31)^2+(-31k^2-60k+1)^2+(-k^2-62k-30)^2+$$
$$(36^2+7^2+17^2)*(k^2+k+1)^2$$
$$= (3^2+4^2+19^2+27^2+34^2+35^2)*(k^2+k+1)^2$$

For k=2 we get the below mentioned new numerical solution:

$$(85,158,243,252,49,119)^2 = (21,28,133,189,238,245)^2$$

---

**Degree three, n=3,**

First method:

$$a^3+b^3+c^3+d^3+e^3+f^3 = p^3+q^3+r^3+s^3+t^3+u^3 \text{ ---- (3)}$$

Let (a,b,c,d,e,f)=[(Ax+1),(Bx+1),(Cx+1),(Dx+1),(Ex+1),(Fx+1)]

and

(p,q,r,s,t,u)=[( Px+1),(Qx+1),(Rx+1),(Sx+1),(Tx+1),(Ux+1)]

Let $(A,B,C,D,E,F)^3 = (P,Q,R,S,T,U)^3 -----(4)$ be known solution

We have numerical solution, $(1,2,4,8,9,12)^3 = (3,5,6,7,10,11)^3$

Substituting values of (a,b,c,d,e,f) &( p,q,r,s,t,u) in equation (3) above we get after simplification the below mentioned condition,

$$x = -\frac{[(A+B+C+D+E+F)-(P+Q+R+S+T+U)]}{W} \quad ----------- (5)$$

Where, W= $(P^2+Q^2+R^2+S^2+T^2+U^2)-(A^2+B^2+C^2+D^2+E^2+F^2)$

After substituting values of $(A,B,C,D,E,F)$ *and* $(P,Q,R,S,T,U)$

in equation (5) we get the solution x = - (1/5)

Hence substituting values of (A, B, C, D, E, F) and $(P,Q,R,S,T,U)$ in equation (3) above we get for x= (-1/5),

$$[(x+1),(2x+1),(4x+1),(8x+1),(9x+1),(12x+1)]^3 =$$
$$[(3x+1),(5x+1),(6x+1),(7x+1),(10x+1),(11x+1)]^3$$

Substituting x= (-1/5) we get,

$$(-4,-3,-1,3,4,7)^3 = (-2,0,1,2,5,6)^3$$

Similarly for different numerical solutions for equation (4) we can get another value of 'x' from equation (5) and thus more numerical solutions.

---

Second Method:

$$a^3+b^3+c^3+d^3+e^3+f^3 = p^3+q^3+r^3+s^3+t^3+u^3 ----(3)$$

Let (a,b,c,d,e,f)=[(Ax+1),(Bx-1),(Cx+1),(Dx-1),(Ex+1),(Fx-1)]

&

(p,q,r,s,tu)= [( Ax-1),(Bx+1),(Cx-1),(Dx+1),(Ex-1),(Fx+1)]

Substituting values of (a,b,c,d,e,f) &( p,q,r,s,t,u) in equation (3) above we get after simplyfication the below mentioned condition,

$(A^2 + C^2 + E^2) = (B^2 + D^2 + F^2)$ ----------- (4)

Equation (4) has the numerical solution, (A,C,E)=(2,10,21) & (B,D,F)=(5,6,22)

Hence substituting values of (A, B, C, D, E, F) we get the parametric form given below,

$$[(2x+1),(5x-1),(10x+1),(6x-1),(21x+1),(22x-1)]^3 =$$

$$[(2x-1),(5x+1),(10x-1),(6x+1),(21x-1),(22x+1)]^3$$

For x=1 we get,

$$(3,4,5,11,22,21)^3 = (1,6,9,20,7,23)^3$$

---

**Degree four, n=4,**

$$a^4 + b^4 + c^4 + d^4 + e^4 + f^4 = p^4 + q^4 + r^4 + s^4 + t^4 + u^4 \text{ ----(5)}$$

We have numerical solution given below,

$$(16,480,496,532,798,1330)^4 = (342,336,224,560,9501292)^4$$

In the above let a=16, b=480, c=496. Since (16+480=496) or (a+b=c)

So we have $(a^4 + b^4 + c^4)$ = $a^4 + b^4 + (a+b)^4$=2$(a^2 + ab + b^2)$^2 -----(6)

Let a=480+t and b=16+kt

Hence parameterization $(a^2 + ab + b^2)$ on the right hand side of equation (6) we get

a= $(496k^2 + 960k - 16)/(k^2 + k + 1)$

b= $(-16k^2 - 992k - 480)/(k^2 + k + 1)$

c= $(480k^2 - 32k - 496)/(k^2 + k + 1)$

$$Since\ (a,b,c)^4 = (16,480,496)^4 = (342,336,224,560,9501292)^4 - (532,798,1330)^4$$

After substituting for (a,b,c) we get the below mentioned parameterization,

$$[(496k^2 + 960k - 16)\text{^}4 + (-16k^2 - 992k - 480)\text{^}4 + (480k^2 - 32k - 496)\text{^}4 +$$

$(532^4 + 798^4 + 1330^4) * (k^2 + k + 1)$^4=

$(342^4 + 336^4 + 224^4 + 560^4 + 950^4 + 1292^4)$*$(k^2 + k + 1)$^4

For k=2 we get,

$$(1944,1264,680,1862,2793,4655)^4 = (1197,3325,1176,4522,784,1960)^4$$

---

**Degree five, n=5,**

$$a^5 + b^5 + c^5 + d^5 + e^5 + f^5 = p^5 + q^5 + r^5 + s^5 + t^5 + u^5 \text{ ----(7)}$$

Consider equation $a^5 + b^5 + c^5 + d^5 + e^5 + f^5 = 2(t)^5 - - - - -$ (8)

Let,

a=$Am^2$+Um+3B

b=$Bm^2$-Um+3A

c=$Cm^2$+Vm+3D

d=$Dm^2$-Vm+3C

e=$Em^2$+Wm+3F

f=$Fm^2$-Wm+3E

t=T$m^2$+3T

Substituting the values of (a,b,c,d,e,f) in equation (8) above we get,

$$(Am^2 + um + 3B)^5 + (Bm^2 - um + 3A)^5 + (Cm^2 + um + 3D)^5 +$$

$$(Dm^2 - vm + 3C)^5 + (Em^2 + wm + 3F)^5 + (Fm^2 - wm + 3E)^5 =$$

$$2(T)^5 * (m^2 + 3)^5 \text{ ----------- (9)}$$

We have known solution given below,

$$(A, B, C, D, E, F)^5 = (91,7,-21,119,161,-63)^5 = 2 * (147)^5 = 2 * (T)^5$$

Since in equation (9) we have the parameter ‘m’, so when we substitute values of

((A,B,C,D,E,F) in equation (9), the only unknowns are (u,v,w) , after simplification of equation (9) we get the below conditions,

(U+V+W) = 4(D-F) &

U=2(D-F)

V=2(2D-3B+F)

W=2(-D+3B-2F)

Since (D,B,F) = (119,7,-63) we get (U,V,W) = (364,308,56)

Substituting the values of (A,B,C,D,E,F,U,V,W,T) in equation (9) above we get the below mentioned parameterization ,

$$(91m^2 + 364m + 21)^5 + (7m^2 - 364m + 273)^5 + (-21m^2 + 308m + 357)^5 + (119m^2 - 308m - 63)^5 +$$

$$(161m^2 + 56m - 189)^5 + (-63m^2 - 56m + 483)^5 = 2(147)^5 * (m^2 + 3)^5 \text{ ------------(10)}$$

We also have the below mentioned numerical solution,

$$(91,7,-21,119,161,-63)^5 = 2 * (147)^5 \text{ and}$$
$$(159,-61,127,-29,81,17)^5 = 2 * (147)^5$$

Using the new values (159,-61,127,-29,81,17) in place of (A,B,C,D,E,F) in equation (9)

we get another set of values for (U,V,W) =(-92,284,-376)

Substituting the above values in equation (9) we get another parametric equation given below,

$$(159m^2 - 92m - 183)^5 + (-61m^2 + 92m + 477)^5 + (127m^2 + 284m - 87)^5 +$$

$$(-29m^2 - 284m + 381)^5 + (81m^2 - 376m + 51)^5 + (17m^2 + 376m + 243)^5 = 2(147)^5 * (m^2 + 3)^5 ----(11)$$

Since the right hand sides of equations (10) & (11) are equal, hence we can equate their left hand sides and we get the below mentioned (5-6-6) equation for degree n=5.

$$(91m^2 + 364m + 21)^5 + (7m^2 - 364m + 273)^5 + (-21m^2 + 308m + 357)^5 + (119m^2 - 308m - 63)^5 +$$

$$(161m^2+56m-189)^5+(-63m^2-56m+483)^5 =$$

$$(159m^2-92m-183)^5+(-61m^2+92m+477)^5+$$

$$(127m^2+284m-87)^5+(-29m^2-284m+381)^5+(81m^2-376m+51)^5+$$

$$(17m^2+376m+243)^5$$

For m=2 we get,

$$(1113,377,889,303,567,119)^5=(269,417,989,203,427,1063)^5$$

---

**Degree six, n=6,**

There are parameter solutions to,

$A_1^6+A_2^6+A_3^6+A_4^6+A_5^6+A_6^6 = B_1^6+B_2^6+B_3^6+B_4^6+B_5^6+B_6^6$ ------ (2)

Let,

$A_1 = a_1a+b_1b-c_1$
$A_2 = a_3a+c_3$
$A_3 = a_4a+c_4$
$A_4 = a_1a-b_2b+c_2$
$A_5 = a_2a-b_1b+c_2$
$A_6 = a_2a+b_2b+c_2$

$B_1 = a_1a+b_1b+c_1$
$B_2 = a_3a-c_3$
$B_3 = a_4a-c_4$
$B_4 = a_1a-b_2b-c_2$
$B_5 = a_2a-b_1b-c_2$
$B_6 = a_2a+b_2b-c_2$

$a_2=ma_1, a_3=na_1, a_4=pa_1, c_1=qc_2, c_3=rc_2, c_4=tc_2$

As for this equation (2), the factorization is done as follows.

(m,n,p,q,r,t) = (3, 1, -1, -2, -6, 3)

$240a_1ac_2(19c_2^2+b_2^2b^2+2a_1b_2ba+4a^2a_1^2)*(21c_2^2+b_2^2b^2+2a_1b_2ba+6a^2a_1^2)$

(m,n,p,q,r,t) = (3, 1, -1, -1, -5, 3)

$240a_1ac_2(14c_2^2+6a_1^2a^2+2b_2ba_1a+b_2^2b^2)*(-12c_2^2+4a_1^2a^2+2b_2ba_1a+b_2^2b^2)$

(m,n,p,q,r,t) = (3, 1, -1, 2, -2, 3)

$240a_1ac_2(5c_2^2+b_2^2b^2+2a_1ab_2b+6a_1^2a^2)*(-3c_2^2+b_2^2b^2+2a_1ab_2b+4a_1^2a^2)$

(m,n,p,q,r,t ) = (3, 1, -1, 3, -3, 1)

$240a_1ac_2(6c_2^2+b_2^2b^2+2a_1ab_2b+6a^2a_1^2)*(-4c_2^2+b_2^2b^2+2a_1ab_2b+4a^2a_1^2)$

(m,n,p,q,r,t) = (3, 2, -4, 3, 2, 2)

$-240a_1ac_2(5c_2^2+15a^2a_1^2+2ba_1ab_2+b_2^2b^2)*(3c_2^2+5a^2a_1^2-2ba_1ab_2-b_2^2b^2)$

(m,n,p,q,r,t) = (3, 3, -5, 3, 2, 2)

$-240a_1ac_2(3c_2^2-b_2^2b^22a_1ab_2b+12a_1^2a^2)*(5c_2^2+b_2^2b^2+2a_1ab_2b+22a_1^2a^2)$

(m,n,p,q,r,t) = (3, 4, -6, 3, 2, 2)

$-240a_1ac_2(3c_2^2+21a_1^2a^2-2b_2a_1ba-b_2^2b^2)*(5c_2^2+31a_1^2a^2+2b_2a_1ba+b_2^2b^2)$

Take $c_2=1$

Case 1. $4a^2a_1^2+2a_1b_2ba+b_2^2b^2=19$.................................. (3)

We can find infinitely many rational solutions of (3), then we obtain infinitely parameter solutions of (2).

Take $x=2a_1a$, $y=b_2b$

Then (3) becomes to $x^2+xy+y^2=19$...................................... (4)

(x,y) = (3,2) is a solution of (4).

We obtain parameter solution of (3) by using $(a,b)=(3/(2a_1),2/b_2)$. And ‘k’ is parameter.

$A_1=7k^2b_2^2-4a_1^2$
$A_2=-9k^2b_2^2-32kb_2a_1-68a_1^2$
$A_3=3k^2b_2^2+20kb_2a_1+44a_1^2$
$A_4=15k^2b_2^2+20kb_2a_1-28a_1^2$
$A_5=11k^2b_2^2-20kb_2a_1-52a_1^2$
$A_6= k^2b_2^2-44kb_2a_1-36a_1^2$

$B_1=-k^2b_2^2-16kb_2a_1-36a_1^2$
$B_2=15k^2b_2^2+16kb_2a_1+28a_1^2$
$B_3=-9k^2b_2^2-4kb_2a_1-4a_1^2$
$B_4=11k^2b_2^2+12kb_2a_1-44a_1^2$
$B_5=7k^2b_2^2-28kb_2a_1-68a_1^2$

$B_6=-3k^2b_2^2-52kb_2a_1-52a_1^2$

Example,

$(a_1,b_{2,}) = (1,1)$ and $(k=1)$

a) $3^6+ 109^6+ 67^6+ 7^6+ 61^6+ 79^6 =$
$53^6+ 59^6+ 17^6+ 21^6+ 89^6+ 107^6$

b) $59^6+ 245^6+ 131^6+ 167^6+ 13^6+ 159^6 =$
$93^6+ 211^6+ 97^6+ 91^6+ 89^6+ 235^6$

c) $27^6+ 85^6+ 43^6+ 73^6+ 11^6+ 49^6 =$
$29^6+ 83^6+ 41^6+ 45^6+ 17^6+ 77^6$

---

**Degree Seven, n=7,**

Let,

$$x_1^k + x_2^k + x_3^k + x_4^k + x_5^k + x_6^k = y_1^k + y_2^k + y_3^k + y_4^k + y_5^k + y_6^k \qquad \text{------- (1)}$$

We show that,

for k=1, 3, 5, 7 eqn. (1) has infinitely many integer solutions.

We used the below identity and Theorem.

Identity (Tito Piezas) (Ref.no.10).

$$(c + bp + aq)^k + (c - bp - aq)^k + (d + ap - bq)^k + (d - ap + bq)^k = (c + bp - aq)^k + (c - bp + aq)^k + (d + ap + bq)^k + (d - ap - bq)^k,$$

for k=2,4,6,
where $(4q^2 - p^2)b^2 + (4p^2 - q^2)a^2 = 15c^2$
$and\ (4p^2 - q^2)b^2 + (4q^2 - p^2)a^2 = 15d^2.$

Theorem,

If $a_1$ ^k+$a_2$ ^k+..........+ $a_m$^k =$b_1$ ^k+ $b_2$^k+............+$b_m$ ^k, for k=2,4,...2n, then(using t) we get,

$$(t+a_1)^k+(t+a_2)^k+\cdots+(t+a_m)^k,(t-a_1)^k+(t-a_2)^k+\cdots+(t-a_m)^k=$$

$$(t+b_1)^k+(t+b_2)^k+\cdots+(t+b_m)^k,(t-b_1)^k+(t-b_2)^k+\cdots+(t-b_m)^k,$$

for k=1,2,3,...2n+1, where t is arbitrary integer.

where $(4q^2-p^2)b^2+(4p^2-q^2)a^2=15c^2$

$$and\ (4p^2-q^2)b^2+(4q^2-p^2)a^2=15d^2$$

Above for k=1,3,5,7 has infinite many solutions,

| | |
|---|---|
| $x_1$ =-c+d+ap-bq | $y_1$ =-c+d+ ap+bq |
| $x_2$=-c+d-ap+bq | $y_2$ =-c+d-ap-bq |
| $x_3$ =-2c-bp-aq | $y_3$ =-2c-bp+aq |
| $x_4$=-2c+bp+aq | $y_4$ =-2c+bp-aq |
| $x_5$ =-c-d-ap+bq | $y_5$ =-c-d-ap-bq |
| $x_6$ =-c-d+ap-bq | $y_6$ =-c-d+ ap+bq |

We have,

$$x_1^k+x_2^k+x_3^k+x_4^k+x_5^k+x_6^k=$$
$$y_1^k+y_2^k+y_3^k+y_4^k+y_5^k+y_6^k \qquad \text{------- (1)}$$

Let [$a_1$, $a_2$,..., $a_m$] = [$b_1 b_2$,..., $b_m$] (k=1,2,...,n) denote

$a_1$ ^k+$a_2$ ^k+..........+ $a_m$^k =$b_1$ ^k+ $b_2$^k+............+$b_m$ ^k (k=1,2,...,n).

First we apply Tito Piezas's identity to above theorem, we obtain (using t)

$$[t+c+bp+aq,t+c-bp-aq,t+d+ap-bq,t+d-ap+bq,t-c-bp-aq,t-c+bp+aq,t-d-ap+bq,t-d+ap-bq]^k=$$

$$[t+c+bp-aq,t+c-bp+aq,t+d+ap+bq,t+d-ap-bq,t-c-bp+aq,t-c+bp-aq,t-d-ap-bq,t-d+ap+bq]^k,$$

where $(4q^2-p^2)b^2+(4p^2-q^2)a^2=15c^2$

$$and\ (4p^2-q^2)b^2+(4q^2-p^2)a^2=15d^2$$

for k=1,3,5,7,

Here, set t = (-c), we can reduce four terms and obtain

$$[-c+d+ap-bq,-c+d-ap+bq,-2c-bp-aq,-2c+bp+aq,-c-d-ap+bq,-c-d+ap-bq]^k =$$

$$[-c+d+ap+bq,-c+d-ap-bq,-2c-bp+aq,-2c+bp-aq,-c-d-ap-bq,-c-d+ap+bq]^k,$$

for k=1, 3, 5, 7
where $(4q^2-p^2)b^2+(4p^2-q^2)a^2=15c^2$
$and\ (4p^2-q^2)b^2+(4q^2-p^2)a^2=15d^2$

By transform simultaneous equation $\{(4q^2-p^2)b^2+(4p^2-q^2)a^2=15c^2$ ,
$(4p^2-q^2)b^2+(4q^2-p^2)a^2=15d^2\}$ $\quad to\ an\ elliptic\ curve,$
we can prove the infinity of solutions to this equation.
Since Tito Piezas's identity has infinitely many integer solutions, therefore this identity also has infinitely many integer solutions.

Example,
(p, q) = (3,2)

$$[-c+d+3a-2b,-c+d-3a+2b,-2c-3b-2a,-2c+3b+2a,-c-d-3a+2b,-c-d+3a-2b]^k =$$

$$[-c+d+3a+2b,-c+d-3a-2b,-2c-3b+2a,-2c+3b-2a,-c-d-3a-2b,-c-d+3a+2b]^k$$

for k=1, 3, 5, 7
where, $7b^2+32a^2=15c^2$
$and$ $32b^2+7a^2=15d^2.$

Numerical example,

$\{a,b,c,d\}=\{1,13,9,19\}$:
$$[-13,33,-59,23,-5,-51]^k = [39,-19,-55,19,-57,1]^k$$
$$\{a,b,c,d\}=\{466,607,797,942\}:$$
$$[329,-39,-4347,1159,-1923,-1555]^k = [2757,-2467,-2483,-705,-4351,873]^k$$
$$\{a,b,c,d\}=\{607,466,942,797\}:$$
$$[372,-517,-2248,364,-1314,-425]^k = [1304,-1449,-1034,-850,-2246,507]^k$$

(p, q)= (4,1)

$$[-c+d+4a-b,-c+d-4a+b,-2c-4b-a,-2c+4b+a,-c-d-4a+b,-c-d+4a-b]^k =$$

$$[-c+d+4a+b, -c+d-4a-b, -2c-4b+a, -2c+4b-a, -c-d-4a-b, -c-d+4a+b]^k$$

for k=1, 3, 5, 7

where, $(-12b^2+63a^2)=15c^2$ and $63b^2-12a^2=15d^2$.

One Solution is,

$$\{a,b,c,d\}=\{89,82,167,148\}:$$

$$x_1^k+x_2^k+x_3^k+x_4^k+x_5^k+x_6^k = y_1^k+y_2^k+y_3^k+y_4^k+y_5^k+y_6^k$$

For k=7, we get:

$$[255, 457, 573, 83, 95, 753]^7 = [419, 293, 751, 589, 41, 123]^7$$

Second solution is,

$$\{a,b,c,d\}=\{82,89,148,167\}:$$

$$[129, 199, 285, 71, 11, 366]^7 = [218, 110, 367, 277, 38, 51]^7$$

---

**Degree eight, n=8,**

$$a^8+b^8+c^8+d^8+e^8+f^8 = p^8+q^8+r^8+s^8+t^8+u^8$$

There are many solutions for $A_1^8+A_2^8+A_3^8+A_4^8+A_5^8+A_6^8 = B_1^8+B_2^8+B_3^8+B_4^8+B_5^8+B_6^8$.

We show that there are infinitely many solutions for 8.6.6 by Sinha's Theorem.

By Sinha's Theorem,

$A_1^8+A_2^8+A_3^8+A_4^8+A_5^8+A_6^8+A_7^8 = B_1^8+B_2^8+B_3^8+B_4^8+B_5^8+B_6^8+B_7^8$.

But, if we set ($A_1+B_1=0$), we can obtain the solution for (8.6.6) equation.

**Theorem**

There are infinitely many solutions for $A_1^8+A_2^8+A_3^8+A_4^8+A_5^8+A_6^8 = B_1^8+B_2^8+B_3^8+B_4^8+B_5^8+B_6^8$.

1. Solving for $a_1^2 + a_2^2 + a_3^2 = b_1^2 + b_2^2 + b_3^2$.

Set $a_1=a(x)+s_1$, $a_2=b(x)+s_2$, $a_3=b(x)+s_2-3(ax+s_1)$,

$b_1=a(x)-s_2$, $b_2=b(x)-s_1$, $b_3=(b-3a)x-s_2+3s_1$......................(1)

Take $s_1=5a-3b$,$s_2=19a-5b$ then

$a_1^2 + a_2^2 + a_3^2 - (b_1^2 + b_2^2 + b_3^2) = 0$

2. Solving for $a_1^4 + a_2^4 + a_3^4 = b_1^4 + b_2^4 + b_3^4$.

$a_1^4 + a_2^4 + a_3^4 -( b_1^4 + b_2^4 + b_3^4) = -32x(a+b)(3a-b)*f$

$f= (8x^2+21x-275)*a^2+(-5x^2-24x+170)*ab+(3x-3)b^2$

We must find rational value (a,b) for above equation.
Discriminant,

$25x^4+144x^3-1280x^2-4608x+25600 = y^2$............................ (2)

So, we must find rational numbers x, y.

U=x and V=y

$V^2 = 25U^4+144U^3-1280U^2-4608U+25600$................. (3)

Using elliptic curve theory, and using 'APECS' program by Ian Connell,
Weierstrass form is,

$Y^2=X^3+X^2-920X+10404$............................ (4)

U = (200X-3248)/(5Y+9X-162)

$V = (25344Y-194880X^2+3550080X-23468800+4000X^3)/(5Y+9X162)^2$.........(5)

$X = (5V+800-72U-7U^2)/U^2$

$Y = (200V+32000-4320U-800U^2-9VU+45U^3)/U^3$....................... (6)

Point P = (0,160) solution for (3).
Rational point Q(X, Y) on the curve (4) corresponding to the values U=0, V=160 is
X=406/25, Y=-396/125.

So, we get the relation of the curve (3) and the curve (4).

Point , (P) =(0,160) on the curve (3), Point, (Q) =(406/25, -396/125) on the curve (4)

We obtain 2Q = (4939/25, -344112/125) on the curve (4) using APECS.

As this point on the curve (4) does not have integer coordinates,
there are infinitely many rational points on the curve (4) by Nagell-Lutz theorem.

Point, ( 2P ) = (-200/67, 725280/4489) is given by 2Q using (5).

We can obtain infinitely many integer solutions for (2) by applying the group law.

By Sinha's Theorem,

$$A_1^8+ A_2^8+ A_3^8+ A_4^8+ A_5^8+ A_6^8+ A_7^8 = B_1^8+ B_2^8+ B_3^8+ B_4^8+ B_5^8+ B_6^8+ B_7^8 \quad \text{---- (1)}$$

$A_1=2a_1$
$A_2=2a_2$
$A_3=b_1+b_2+b_3$
$A_4=2a_3$
$A_5=b_1-b_2+b_3$
$A_6=-b_1+b_2+b_3$
$A_7=b_1+b_2-b_3$

$B_1=a_1-a_2+a_3$
$B_2=-a_1+a_2+a_3$
$B_3=2b_3$
$B_4=a_1+a_2+a_3$
$B_5=2b_1$
$B_6=2b_2$
$B_7=a_1+a_2-a_3$

By Sinha's Theorem, substitute (a, b, x) to (1), then
we obtain infinitely many solutions of (1).

Example,

(x, a, b ) : $(8x^2+21x-275)a^2+(-5x^2-24x+170)*ab+(3x-3)b^2=0$

Since, $a_1=ax+s_1$, $a_2=bx+s_2$, $a_3=bx+s_2-3(ax+s_1)$, $b_1=ax-s_2$, $b_2=bx-s_3$,
$b_3 = (b-3a)x-s_2+3s_1$……………………………….....(2)

and, $s_1$=5a-3b, $s_2$=19a-5b then

For, (x, a, b) =(1 , 47 ,82) we get after substituting values of ($a_1$ ,$a_2$ ,$a_3$ ,$b_1$ ,$b_2$ ,$b_3$ ) in equation (2) and equation (1) above,

$565^8+ 459^8+ 457^8+ 552^8+ 23^8+ 116^8=$
$493^8+ 575^8+ 529^8+ 436^8+ 93^8+ 72^8$

(x, a, b ) = (-6, 21, 113)

$211^8+ 155^8+ 59^8+ 44^8+ 165^8+ 54^8=$
$31^8+ 209^8+ 121^8+ 10^8+ 111^8+ 180^8$

(x, a, b) = (6, 15, 139)

$106^8+ 203^8+ 295^8+ 91^8+ 78^8+ 216^8=$
$232^8+ 13^8+ 169^8+ 125^8+ 294^8+ 126^8$

(x, a, b) = (-14, 5, 9)

$19^8+ 27^8+ 35^8+ 4^8+ 3^8+ 34^8=$
$17^8+ 7^8+ 1^8+ 30^8+ 31^8+ 36^8$

(x, a, b) = (-14,3 , -37)

$190^8+ 111^8+ 127^8+ 13^8+ 182^8+ 84^8=$
$148^8+ 195^8+ 169^8+ 71^8+ 98^8+ 42^8$

---------------------------------------------------------------------------------------------

**<u>Degree nine , n=9,</u>**

$$x_1^9 + x_2^9 + x_3^9 + x_4^9 + x_5^9 + x_6^9 = y_1^9 + y_2^9 + y_3^9 + y_4^9 + y_5^9 + y_6^9$$

Andrew Bremner and J. Delorme (Ref. no. 1) showed that,

$$x_1^k + x_2^k + x_3^k + x_4^k + x_5^k + x_6^k = y_1^k + y_2^k + y_3^k + y_4^k + y_5^k + y_6^k$$

for k=1, 2, 3, 9 has infinitely many integer solutions.([1])

By differrent way, Tito Piezas [10] showed that above equation has infinitely many solutions.

A.Choudhry (Ref. no. 11) showed that $x_1^k + x_2^k + x_3^k = y_1^k + y_2^k + y_3^k$

for k=1,2,6 has infinitely many integer solutions ([3])

We show that there are infinity many solutions for above (k.6.6) equation.

$$x_1^k + x_2^k + x_3^k + x_4^k + x_5^k + x_6^k = y_1^k + y_2^k + y_3^k + y_4^k + y_5^k + y_6^k \text{------- (2)}$$

for k=1,2,3,9 has infinitely many solutions**.**

Set variables as following [equation (3)],

$x_1$ = 2(a+b)m+(a-b+t)n+w | $y_1$ =2(a+b)m+(a-b+t)n-w

$x_2$= -2am+(a+b+t)n+w | $y_2$ =-2am+(a+b+t)n-w

$x_3$ = -2bm-(a+b-t)n+w | $y_3$ =-2bm-(a+b-t)n-w

$x_4$= -2(a+b)m+(a-b+t)n-w | $y_4$ =-2(a+b)m+(a-b+t)n+w

$x_5$ = 2am+(a+b+t)n-w | $y_5$ =2am+(a+b+t)n+w

$x_6$ = 2bm-(a+b-t)n-w | $y_6$ = 2bm-(a+b-t)n+w

By using Ajai Choudhry's [Ref. no.11], equation (2) is always equals to zero for degree k=1, 2, 3.

And for k=9 we have :

$$x_1^9 + x_2^9 + x_3^9 + x_4^9 + x_5^9 + x_6^9 = y_1^9 + y_2^9 + y_3^9 + y_4^9 + y_5^9 + y_6^9$$

Substituting values of above from equation (3) and after simplification we get,

Expression=18432*a*b*m* $n^2$(m-n)(m+n)(a+b)(a-b+3t)

$$* (-5m^2b^3 - n^2b^3 + 7m^2b^2t - 4m^2b^2a + 7n^2b^2t + 2n^2b^2a + 7am^2tb + 4m^2ba^2 - 14n^2bt^2 - 2n^2ba^2 - 7an^2bt + 5m^2a^3 + 7m^2a^2t + 7n^2a^2t + 14n^2at^2 + n^2a^3 + 14n^2t^3)$$
$$* (2m^2b^2 + 2m^2ba + 2m^2a^2 + 3n^2a^2 + 14n^2ta + 21n^2t^2 - 4bn^2a - 14bn^2t + 3b^2n^2)$$

So, we have to find the rational solution (m, n) of:

$(-14bt^2 - 7bat + a^3 + 2ab^2 + 14at^2 - 2a^2b - b^3 + 7a^2t + 7b^2t + 14t^3)n^2 +$
$(-4ab^2 - 5b^3 + 4a^2b + 5a^3 + 7b^2t + 7a^2t + 7bat)m^2 = 0$........... (3)

So that there are rational solution,

$$-(-14bt^2 - 7bat + a^3 + 2ab^2 + 14at^2 - 2a^2b - b^3 + 7a^2t + 7b^2t + 14t^3) * (-4ab^2 - 5b^3 + 4a^2b + 5a^3 + 7b^2t + 7a^2t + 7bat)$$

must be square number ( $s^2$), then we have to find rational solution (a, b, t, s) of

$$s^2 = (-98a^2 - 98ab - 98b^2)t^4 + (56ab^2 - 56a^2b + 168b^3 - 168a^3)t^3 + (63a^2b^2 - 119a^4 + 14ab^3 - 119b^4 + 14a^3b)t^2 + (14a^3b^2 - 14a^2b^3 - 42a^5 + 42b^5 - 14ab^4 + 14a^4b)t + 6ab^5 + 6a^5b - 5a^6 + 2a^4b^2 - 5b^6 - 6a^3b^3 + 2a^2b^4 \ldots..(4)$$

By computer search, we found a solution (a, b, t)=(3, 4, 27/41).

Substitute (a, b, t)=(3, 4, 27/41) to (3),then (3) becomes to 160/68921(139n-164m)(139n+164m)=0.

So, we get (m, n)=(139,164).

Substitute (a, b, t, m, n) = (3, 4, 27/41,139,164) to (2),
then we get following solution,

[x1, x2, x3, x4, x5, x6] = [1025, 291, -996, -1081, 965, -44]
[y1, y2, y3, y4, y5, y6] = [865, 131,-1156, -921, 1125, 116].

Next, substitute ,

(a, b) = (3, 4) to (4),then we get a quartic equation,

$s^2$= -3626*t^4+6888*t^3-26831*t^2+24570*t-3029............. (5)

Using elliptic curve theory, transform (5) to minimal Weierstrass form (6).

$V^2$ + UV + V = $U^3$ -7166374 -22875861928............. (6)

We get a point P (U, V) = (1026337/64, -1026359837/512).

As this point on the curve (6) does not have integer coordinates,
there are infinitely many rational points on the
curve (6) by Nagell-Lutz theorem.

By using point 2P= (t, s)=(3181201/12876603, 6408411316637440/165806904819609),
we obtain a new solution.

[x1, x2, x3, x4, x5, x6] =
[15677071397, 40208111671,-63297775068,-26458358421, 63560861593,-33396207172]
[y1, y2, y3, y4, y5, y6] =
[19383367397, 43914407671, -59591479068,-30164654421, 59854565593,-37102503172]

Example,
[ a, b, t] = [1, 3, 6/5],

$$[\ 18,13,14,23,13,1]^9 = [\ 5,10,15,21,22,9]^9$$

[a,b,t] = [4, 9, 13/3],

$$[\ 453,122\ ,331,431\ ,150,281]^9 = [\ 429,98\ ,307\ ,455\ ,174,305]^9$$

---

**<u>Table (A) :</u>**

$(\ six)\ n^{th}\ powers\ equal\ to\ another\ (six)\ n^{th}\ powers$ :

$$[a^n + b^n + c^n + d^n + e^n + f^n] = [p^n + q^n + r^n + s^n + t^n + u^n]$$

Numerical solutions: For degree's n=2, 3, 4, 5, 6, 7, 8 & 9

| n | a | b | c | d | e | f | | p | q | r | s | t | u |
|---|---|---|---|---|---|---|---|---|---|---|---|---|---|
| 2 | 1 | 7 | 17 | 30 | 31 | 36 | = | 3 | 4 | 19 | 27 | 34 | 35 |
| 3 | 11 | 22 | 4 | 3 | 21 | 5 | = | 20 | 7 | 6 | 23 | 9 | 1 |
| 4 | 16 | 480 | 496 | 532 | 798 | 1330 | = | 224 | 342 | 336 | 560 | 950 | 1292 |
| 5 | 87 | 233 | 264 | 396 | 496 | 540 | = | 90 | 206 | 309 | 366 | 522 | 523 |
| 6 | 61 | 3 | 109 | 67 | 7 | 79 | = | 21 | 17 | 53 | 59 | 89 | 107 |
| 7 | 129 | 199 | 285 | 71 | 11 | 366 | = | 218 | 110 | 367 | 277 | 38 | 51 |
| 8 | 3 | 6 | 8 | 10 | 15 | 23 | = | 5 | 9 | 12 | 9 | 20 | 22 |
| 9 | 1 | 13 | 14 | 13 | 18 | 23 | = | 5 | 9 | 10 | 15 | 21 | 22 |